\documentclass[12pt]{amsart}
\usepackage{amsmath}
\usepackage{amssymb}
\usepackage{epsfig}
\usepackage{graphicx}
\numberwithin{equation}{section}

\input xy
\xyoption{all}

\calclayout
\allowdisplaybreaks[3]

\theoremstyle{plain}
\newtheorem{prop}{Proposition}

\newtheorem{theo}[prop]{Theorem}
\newtheorem{coro}[prop]{Corollary}

\newtheorem{lemm}[prop]{Lemma}

\theoremstyle{definition}

\newtheorem{conj}[prop]{Conjecture}

\newtheorem{rema}[prop]{Remark}

\newtheorem{exam}[prop]{Example}

\def\ocM{\overline {\mathcal M}}
\def\oM{\overline {M}}
\def\oS0{\overline{S_0}}

\def\iHom{{\mathcal H}\!{\it om}}
\def\Hilb{\mathrm{Hilb}}

\def\iExt{{\mathcal E}\!{\it xt}}

\def\mo{{\mathfrak o}}

\def\mp{{\mathfrak p}}

\def\cC{{\mathcal C}}

\def\cE{{\mathcal E}}

\def\cO{{\mathcal O}}

\def\cK{{\mathcal K}}
\def\cL{{\mathcal L}}

\def\cI{{\mathcal I}}

\def\cS{{\mathcal S}}
\def\cT{{\mathcal T}}
\def\cW{{\mathcal W}}

\def\cY{{\mathcal Y}}
\def\cW{{\mathcal W}}
\def\cZ{{\mathcal Z}}

\def\ra{\rightarrow}

\def\F{{\mathbb F}}
\def\P{{\mathbb P}}
\def\Q{{\mathbb Q}}

\def\Z{{\mathbb Z}}

\def\bP{{\mathbb P}}
\def\bQ{{\mathbb Q}}
\def\bH{{\mathbb H}}
\def\bZ{{\mathbb Z}}
\def\bR{{\mathbb R}}

\def\bC{{\mathbb C}}

\def\Pic{{\rm Pic}}

\def\ra{\rightarrow}

\def\Spec{{\rm Spec}}
\def\Specf{{\rm Spf}}

\def\Z{{\mathbb Z}}

\def\cN{{\mathcal N}}
\def\cI{{\mathcal I}}

\def\Pic{{\rm Pic}}

\makeatother
\makeatletter

\author{Fedor Bogomolov}
\address{Courant Institute\\
                New York University \\
                 New York, NY 10012 \\
                USA }
\email{bogomolo@cims.nyu.edu}

\author{Brendan Hassett}
\address{Department of Mathematics\\
Rice University, MS 136 \\
Houston, Texas  77251-1892 \\
USA}
\email{hassett@rice.edu}

\author{Yuri Tschinkel}
\address{Courant Institute\\
                New York University \\
                New York, NY 10012 \\
                USA }
\email{tschinkel@cims.nyu.edu}

\title[Rational curves]{Constructing rational curves on K3 surfaces}

\begin{document}
\date{\today}

\begin{abstract}
We develop a mixed-characteristic version of the Mori-Mukai technique
for producing rational curves on K3 surfaces.  
We reduce modulo $p$, produce rational
curves on the resulting K3 surface over a finite field, and lift to
characteristic zero.
As an application, we prove that all complex K3 surfaces with Picard group generated by a class of
degree two have an infinite
number of rational curves.  
\end{abstract}

\maketitle
\tableofcontents

\section{Introduction}
\label{sect:introduction}

Let $K$ be an algebraically closed field and 
$S$ a K3 surface defined over $K$. 
It is known that $S$ contains rational curves--see Mori-Mukai \cite{mori-mukai},
as well as Theorem~\ref{theo:MM} and Corollary~\ref{coro:MM} below.
In fact, an extension of the argument in \cite{mori-mukai} shows that the {\em general} 
K3 surface of given degree has infinitely many rational curves; we
sketch this below in Theorem~\ref{theo:infinite} (cf. \cite{XiChen}).
The idea is to specialize the K3 surface $S$ to a K3 surface $S_0$ with Picard
group of rank 2, where some multiple of the polarization can be expressed as a 
sum of linearly independent classes of smooth rational curves. The union of these
rational curves deforms to an irreducible rational curve on $S$. This idea
applies to K3 surfaces parametrized by points 
outside a countable union of subvarieties of the moduli space. 
In particular, {\em a priori} it doesn't apply to
K3 surfaces over countable fields, 
such as ${\overline \F_p}$ and ${\overline \Q}$.  Of course, there are also other techniques
proving density of rational curves on special K3 surfaces, e.g., certain Kummer surfaces 
\cite{mori-mukai}, surfaces with infinite automorphisms
\cite[proof of Thm. 4.10]{bt}, or with elliptic fibrations
(see Remark~\ref{rema:ellip} below).
These K3 surfaces have Picard rank $\ge 2$, 
and all except finitely many lattices in rank $\ge 3$ correspond to K3 surfaces
with infinite automorphisms or elliptic fibrations \cite{Nik,Vin}.  

Moreover, in \cite{bt-k3} it is proved that, over $k={\overline {\F}_p}$, every
algebraic point on a Kummer K3 surface lies on an irreducible rational curve. 
The proof of this result uses the Frobenius endomorphism on the covering 
abelian surface.

\begin{theo}
\label{thm:main}
Let $S$ be a K3 surface over an algebraically closed field of
characteristic zero with $\Pic(S)=\Z$, generated by a 
divisor of degree two.  
Then $S$ contains infinitely many rational curves. 
\end{theo}

The motivation for our argument comes from a result of Bogomolov and Mumford
\cite{mori-mukai}:  Let $(S,f)$ be a general K3 surface of degree $2g-2$.
We can degenerate $S$ to a Kummer K3 surface $(S_0,f)$, which has 
infinitely many rational curves.  Indeed, we can produce examples where
there are infinitely many (reducible) rational curves in $|Nf|,N\ge 1$,
consisting of unions of smooth components meeting transversally.
A deformation argument shows that these deform to infinitely many
(irreducible) rational curves in nearby fibers.  However, on subsequent
specializations, distinct rational curves might collapse onto each
other.  If there were an infinite number of such collisions, the
specialized K3 surfaces might only have a finite number of rational
curves.  

Here we emulate the argument in \cite{mori-mukai} in
mixed characteristic.  K3 surfaces over finite fields play the
r\^ole of the Kummer surface;  the `general' K3 surface is a 
K3 surface over a number field with Picard group of rank one.  
The main technical issue is that we cannot assume {\em a priori}
that the rational curves on the reduction $\pmod{\mp}$ have mild
singularities.  Thus we are forced to use more sophisticated
deformation techniques.  

\

\noindent {\bf Acknowledgments:}  The first author was supported 
by NSF grant 0701578;  the second author was supported by NSF
grant 0554491;
the third author was supported by NSF grants 0554280 and 0602333.
We are grateful to Yuri Zarhin for useful discussions and comments.  

\section{Guiding questions and examples}
The following is well-known but hard to trace in the literature:

\begin{conj}[Main conjecture] \label{conj:rational}
Let $K$ be an algebraically closed field of arbitrary characteristic and
$S$ a projective K3 surface over $K$.  
There exist infinitely many rational curves on $S$.
\end{conj}

In characteristic zero, we can reduce this to the case of number fields:
\begin{theo}
Assume that for every K3 surface $S_0$ defined over a number field $K_0$,
there are infinitely many rational curves in 
$$\oS0:=S_0\times_{\Spec(K_0)}\Spec(\overline{\bQ}).$$
Then Conjecture~\ref{conj:rational} holds over fields of characteristic zero.
\end{theo}
\begin{proof}
Let $S$ be a K3 surface defined over a field of characteristic zero, which we may
assume is the function field of a variety $B$ defined over a number field $F$.  
Shrinking $B$ as necessary, we obtain a smooth projective morphism
$$\pi:\cS \ra B$$
with generic fiber $S$.  

We claim there exists a point $b\in B({\overline \bQ})$ such that the specialization
map to the fiber $S_b=\pi^{-1}(b)$
$$\Pic(S) \ra \Pic(S_b)$$
is surjective.  The argument is essentially the same as the proof of the main result
of \cite{Ellenberg}.  The only difference is that Ellenberg considers the Galois representation
on the full primitive cohomology of a polarized K3 surface surface, whereas here we
restrict to the representation on the transcendental cohomology of $S$, i.e.,
the orthogonal complement to $\Pic(S)\subset H^2(S,\bZ)$.  

Our assumption is that $S_b$ admits infinitely many rational curves.  We claim each of these
lifts to a rational curve of $S$, perhaps after a generically-finite base-change $\tilde{B}\ra B$.  
Suppose we have a morphism $\phi_b:\bP^1 \ra S_b$, birational onto its image.  The
class $\phi_*[\bP^1]$ remains algebraic in the fibers of $\cS\ra B$.
Thus we can apply a result of Ran \cite[Cor. 3.2 and 3.3]{Ran},
which builds on earlier work of Voisin \cite{Voisin} and Bloch \cite{Bloch}, to conclude that
$\phi_b$ lifts to a morphism $\phi:\bP^1 \ra S$.  
\end{proof}

\begin{rema}
We do not know whether the positive-characteristic case of Conjecture~\ref{conj:rational}
can be reduced to the case of finite fields.
\end{rema}

\begin{exam} \label{exam:Ku}
Here we show that any Kummer K3 surface over an arbitrary algebraically-closed 
field of characteristic $\neq 2$ admits an infinite number of rational curves.

Let $A$ be an abelian surface with Kummer surface $S$:
$$\begin{array}{ccc}
 &   & S \\
 &   & \downarrow \\
A&\ra & A/\pm
\end{array}
$$
Note that $A$ is isogenous
to the Jacobian $J$ of a genus two curve $C$.  (Every abelian surface is isogenous
to a principally-polarized surface, which is either a Jacobian or a product $E_1 \times E_2$ of elliptic
curves.  In the latter case, if we express $E_1$ and $E_2$ as branched coverings on $\bP^1$
at $\{0,\infty, \alpha_1,\beta_1 \}$ and $\{0,\infty,\alpha_2,\beta_2 \}$ with the $\alpha_i$ and
$\beta_i$ distinct, 
then the genus-two
double cover $C\ra \bP^1$ branched at $\{0,\infty,\alpha_1,\beta_1,\alpha_2,\beta_2\}$ 
will work.  Indeed, $E_1$ and $E_2$ are Prym varieties of $C$.)

We choose the embedding $C \hookrightarrow J$ such that a Weierstrass point is mapped
to zero.  Then the images of $n\cdot C$ in $A/\pm$ are distinct rational curves.  Indeed,
multiplication-by-$n$ commutes with $\pm$, and acts on $C$ via the hyperelliptic involution.  
\end{exam}

\begin{rema} \label{rema:ellip}
Elliptic complex K3 surfaces always have infinitely many rational curves:
see \cite[Thm. 1.8]{bt} or \cite[Cor. 8.12, Prop. 9.10, Rem. 9.7]{Ha} and Example~\ref{exam:Ku} above
for the degenerate case where the elliptic surface arises from a Kummer construction.
\end{rema}

\begin{exam}
Let $S$ be a K3 surface over an algebraically-closed field
such that the Picard group is generated over $\bQ$ by 
smooth rational curves $C_1$ and $C_2$ satisfying
\begin{equation} \label{eqn:lattice}
\begin{array}{c|cc}
   & C_1 & C_2 \\
\hline
C_1 & -2 & 6 \\
C_2 & 6 & -2
\end{array}
\end{equation}
and generated over $\bZ$ by $C_1$ and $f=(C_1+C_2)/2$.  Note that
$(S,f)$ can be realized geometrically as the 
double cover of $\bP^2$ branched over a plane sextic curve that
admits a six-tangent conic.  

Surfaces of this type can be defined over $\F_3$ \cite[Ex. 6.1]{EJ}, e.g.,
$$w^2=(y^3-x^2y)^2 + (x^2+y^2+z^2)(2x^3y+ x^3 z+ 2x^2yz+ x^2z^2+2xy^3+2y^4+z^4).$$
The technique of \cite{Ellenberg} can be used to obtain examples over $\bQ$.  
Indeed, the moduli space of lattice-polarized K3 surfaces of type (\ref{eqn:lattice})
is unirational:  The 
sextic plane curves six-tangent
to a fixed conic plane curve $D$ are parametrized by a $\bP^{15}$-bundle over 
$\bP^6$, and these dominate our moduli space.  We can apply Ellenberg's Hilbert irreducibility argument
\cite{Ellenberg} directly to this rational variety.

We do not know how to construct infinitely many rational curves on K3 surfaces of this type, over 
${\overline \bQ}$ or ${\overline \F}_p$.  
\end{exam}

\section{Background results}
A polarized K3 surface $(S,f)$ consists of a K3 surface and
an ample divisor $f$ that is primitive in the Picard group.
Its degree is the positive even integer $f\cdot f$.
Let $\cK_g, g\ge 2$ denote the moduli space (stack) of complex polarized
K3 surfaces of degree $2g-2$, which is
smooth and connected of dimension $19$.

The following result was initially presented by the first author in October 1981 at Mori's
seminar at IAS;  the proof was based on deformation-theoretic ideas developed several years
earlier.  A different argument was presented in \cite{mori-mukai};  Mori and Mukai 
indicate that Mumford also had a proof.  
\begin{theo} \label{theo:MM}  \cite{mori-mukai}
Every complex projective K3 surface contains at least one rational curve.
\end{theo}
For our purposes, it is useful to recall the Mori-Mukai argument for
the existence of rational curves in the generic K3 surface
of degree $2g-2$.

\begin{proof}
We exhibit a K3 surface $S_0$ containing two smooth rational
curves $C_1$ and $C_2$ meeting transversally at $g+1$ points,
such that the class $f=[C_1+C_2]$ is primitive.  
We then deform $C_1 \cup C_2$ to an irreducible rational curve in 
a nearby polarized K3 surface.  

Let $E_1$ and $E_2$ be elliptic curves admitting an isogeny $E_1 \ra E_2$ of
degree $2g+3$ with graph $\Gamma \subset E_1 \times E_2$;
$p\in E_2$
a $2$-torsion point.
Take the associated Kummer surface $S_0$, i.e., the minimal desingularization of
$$(E_1 \times E_2)/\left< \pm 1\right>.$$
Note that $\Gamma$ intersects $E_1 \times p$ transversally in $2g+3$ points, one of which
is $2$-torsion in $E_1 \times E_2$

Take $C_1$ and $C_2$ to be the images of $\Gamma$ and
$E_1 \times p$ in $S_0$, smooth rational curves meeting
transversally in $g+1$ points.
The sublattice of $\Pic(S_0)$ determined by $C_1$ and $C_2$ is:
\begin{equation} \label{K3type}
\begin{array}{c|cc}
   & C_1 & C_2 \\
\hline
C_1 & -2 & g+1 \\
C_2 & g+1 & -2
\end{array}
\end{equation}

Consider deformations of $S_0$ in $\cK_g$,
$$\pi:(\cS,f) \ra  \Delta, \quad \Delta=\{t:|t|<1 \},$$
i.e., deformations for which $f=[C_1]+[C_2]$
remains algebraic.  Since $\cO_{S_0}(C_1+C_2)$ is nef
and big, it has trivial higher cohomology (by Kawamata-Viehweg vanishing)
and $\Gamma(\cO_{S_t}(f))$ is constant in $t$.  
Thus $C_1 \cup C_2$ deforms to divisors in nearby fibers in a smooth family
of relative dimension $g$.  
We can assume that the general fiber of
$\pi$ has Picard group generated by $f$, so that $C_1 \cup C_2$ deforms
to {\em irreducible curves} in the general fiber.  

Consider the rational map
$$\begin{array}{ccc}
\bP(\pi_*\cO_{\cS}(f)) & \stackrel{\mu}{\dashrightarrow} & \oM_g \\
\downarrow & & \\
   \Delta  & & 
\end{array}
$$
assigning to each curve the corresponding point in moduli.  
Note that $\mu$ is regular at $[C_1 \cup C_2]$.  

Let $T$ denote the union of two smooth rational curves at a node, i.e., the
singular conic curve $\{xy=0\}$.  Choose a birational morphism $\phi:T \ra C_1\cup C_2$,
i.e., one that normalizes all but one of the nodes.  
The preimages of the remaining $g$ nodes yield $2g$ smooth points $p_1,\ldots,p_{2g} \in T$,
numbered so that $p_{2i-1}$ and $p_{2i}$ are identified.  
More generally, consider the morphism
$$\iota: \oM_{0,2g} \ra \oM_g$$
which identifies $p_{2i-1}$ and $p_{2i}$ for $i=1,\ldots,g$, so that 
$\iota(T,p_1,\ldots,p_{2g})=C_1 \cup C_2$.  The image of $\iota$ has dimension
$2g-3$ and parametrizes curves with $g$ nodes.  Finally, let $\delta \subset \oM_{0,2g}$
denote the irreducible boundary divisor parametrizing curves with combinatorial type
equal to the type of  $T$.  

Let $Z$ be the closed image of $\mu$, which has dimension $\le g+1$.  Note that
$\mu^{-1}(C_1 \cup C_2)$ is zero-dimensional, in fact, $\iota^{-1}(\mu^{-1}(\delta))$ is zero-dimensional.
Indeed, the generic fibers admit only irreducible curves in $|f|$, and the special
fiber admits only finitely-many rational curves.  Consequently, $\iota(\oM_{0,2g})$ and
$Z$ meet properly and thus in dimension one.  These are stable reductions
of curves appearing in the fibers of $\pi$, which are necessarily of genus zero.
It follows that $\cS$ contains a one-parameter family of genus-zero curves containing
$C_1 \cup C_2$.  Again, since each fiber admits only finitely-many rational curves, 
there exist rational deformations of $C_1 \cup C_2$ in the generic fiber of $\pi$.
\end{proof}

\begin{rema} \label{rema:MM}
This argument requires that the curve in the special K3 surface $S_0$ be nodal.
Without this, we cannot carry out the dimension estimates.  

There exists an (irreducible) eighteen-dimensional family of lattice-polarized K3 surfaces
of type (\ref{K3type}).  {\em A posteriori}, we know that for the generic such K3,
the rational curves $C_1$ and $C_2$ meet transversally.  The argument above is applicable
to any such K3 surface, not just the Kummer surfaces.  

Moreover, even when $C_1$ and $C_2$ fail to meet transversally, $C_1\cup C_2$ is
still the limit of irreducible rational curves.  
\end{rema}

A fairly straightforward specialization argument allows us to deduce that every indecomposable
effective class in a K3 surface contains rational curves.  Any cycle of curves that is a specialization
of a rational curve is a union of rational curves, perhaps with multiplicities.  
We next show that the general K3 surface in $\cK_g$ admits an infinite number of rational curves.

The following theorem has been known to experts;
the first published proof is in \cite{XiChen}.  (Xi Chen attributes
a special case to S. Nakatani.)  
\begin{theo} \label{theo:infinite} \cite{XiChen}
Fix $N\ge 1$.  Then for a generic $(S,f) \in \cK_g$ there exists
an irreducible rational curve in $|Nf|$.  
\end{theo}
\begin{coro} A very general K3 surface of degree $2g-2$ contains
an infinite number of rational curves.
\end{coro}

The proof in \cite{XiChen} involves specializing the K3 surface to a
union of two rational normal scrolls, meeting transversely along an 
elliptic curve.  Xi Chen identifies reducible rational curves on this surface that
can be deformed back to irreducible rational curves on a general
K3 surface.

For our purposes this argument is not sufficiently flexible.  Our technique
entails analyzing the reductions $\pmod{p}$ of a K3 surface defined
over a number field.  We cannot expect these reductions
to be unions of rational normal scrolls.  

We sketch an alternative proof for Theorem~\ref{theo:infinite}, with a view
towards highlighting the geometric ideas behind our main theorem:
\begin{proof}
Let $(S_0,f)$ be a polarized K3 surface with
$$\Pic(S_0)_{\bQ}=\bQ C_1 + \bQ C_2,$$
where $C_1$ and $C_2$ are smooth rational curves such that
$$\begin{array}{c|cc}
   & C_1 & C_2 \\
\hline
C_1 & -2 & N^2(g-1)+2 \\
C_2 & N^2(g-1)+2 & -2
\end{array}.
$$
Furthermore, we assume that 
$$\Pic(S_0)=\bZ C_1 + \bZ C_2 + \bZ f$$
where 
$$Nf=C_1+C_2.$$
The existence of these can be deduced from surjectivity of Torelli.  

In this case, we lack general tools for producing examples where $C_1$
and $C_2$ meet transversally.  This necessitates a modified 
approach to the deformation theory of the union $C_1 \cup C_2$.  

Let $T_0$ again denote the nodal rational curve $\{xy=0 \} \subset \bP^2$,
and choose a birational morphism
$$T_0 \ra C_1 \cup C_2,$$
which maps the node of $T_0$ to some singularity of $C_1\cup C_2$.  
We regard the composition with the inclusion into $S_0$
as a stable map $\phi  \in \ocM_0(S_0,Nf)$.
Again, consider deformations of $S_0$ 
$$\pi:(\cS,f) \ra \Delta,$$
where the generic fiber is a general K3 surface in $\cK_g$,
i.e., one with Picard group generated by $f$.
This time, we consider the relative Kontsevich space of stable maps
$$\psi:\ocM_0(\cS/\Delta,Nf) \ra \Delta,$$
whose fibers are Kontsevich spaces of the fibers of $\pi$.
Recall that a {\em stable map} into a variety $Y$ is a pair $(T,\phi:T \ra Y)$
consisting of a connected nodal curve $T$ and a morphism $\phi$ such
that the dualizing sheaf $\omega_T$ is ample relative to $T$.  

We claim 
that $\ocM_0(\cS/\Delta,Nf)$ is (at least) one-dimensional near $[\phi]$.
Since every fiber of $\pi$ is not uniruled, it follows that
$\psi$ is dominant.  Furthermore, since $C_1$ and $C_2$ do not deform
to algebraic classes in the generic fiber, the resulting genus zero
stable maps have $\bP^1$ as their domain.  They are birational onto
their images, which are therefore irreducible rational curves in $|Nf|$.  
\end{proof}

It remains to prove the claim, which requires some deformation theory, 
explained in the next section.

\section{Deformation results for stable maps}
In this section, we work over a field of arbitrary characteristic.  
Let $Y$ be a smooth projective variety and $\beta$ a curve class on $Y$.

Consider the open substack
$$\ocM^{\circ}_0(Y,\beta)
\subset \ocM_0(Y,\beta)$$
corresponding to maps $\phi:T\ra Y$ that are generic embeddings,
i.e., for each irreducible
component $T_i \subset T$ the restriction of $\phi$ to the generic point of $T_i$ is an
embedding.

This open set dominates the locus $\Xi \subset \Hilb$ consisting of curves $C$
expressible as unions 
$$C=C_1 \cup \ldots \cup C_r \subset Y, \quad \sum_{i=1}^r [C_i]=\beta,$$
of rational curves $C_i$ with each component having multiplicity one;  the induced morphism
is finite-to-one.  Indeed, there are a finite number of connected seminormal curves $T$ through
which the normalization factors
$$C^{\nu}=\coprod_{j=1}^r \bP^1 \ra T \ra C.$$

The substack $\ocM^{\circ}_0(Y,\beta)$ has several nice properties:  First, its objects
have trivial automorphisms, thus the substack is actually a scheme.  (Indeed, it is a finite
cover of $\Xi$, which is quasi-projective.)  
The obstruction theory over $\ocM^{\circ}_{0,0}(Y,\beta)$ takes a particularly simple form:
Given a stable map $\phi:T \ra Y$,
first-order deformations and obstructions are given by 
$$\bH^i( \bR\iHom_{\cO_T}(\Omega_{\phi}^{\bullet},\cO_T)), i=1,2$$
where $\Omega_{\phi}^{\bullet}$ is the complex
$$d\phi^t: \phi^*\Omega^1_Y \ra \Omega^1_T$$
supported in degrees $-1$ and $0$.  
We shall require the following result (cf. \cite[p. 61]{GHS}):

\begin{lemm}
Let $\phi:T \ra Y$ be a stable map to a smooth variety that is unramified at the 
generic point of each irreducible component of $T$.  Then the complex
$$\bR\iHom_{\cO_T}(\Omega_{\phi}^{\bullet},\cO_T)$$
is quasi-isomorphic to $\cN_{\phi}[-1]$, the normal sheaf 
shifted by $-1$.
\end{lemm}
In the special case where the domain $T$ is
smooth, $\cN_{\phi}$ is the cokernel of the differential $d\phi:\cT_T \ra \phi^* \cT_Y$. 
First order deformations of $\phi$ are given by $H^0(\cN_{\phi})$;  obstructions are
given by $H^1(\cN_{\phi})$.  
\begin{proof}
Let $T$ be a nodal projective curve.  Then $\Omega^1_T$ admits a resolution
$$0 \ra \cE_1 \stackrel{f_1}{\ra} \cE_0 \ra \Omega^1_T \ra 0$$
where $\cE_1$ is invertible and $\cE_0$ is locally free of rank two.  Locally,
this takes the following form:  Suppose that $p \in T$ is a node expressed as
$xy=0$ in local \'etale/analytic coordinates.  Then locally
$$\Omega^1_T=(\cO_T dx + \cO_T dy)/\left<ydx+xdy\right>$$
which has a presentation of the specified form.  More intrinsically, we 
may identify $\Omega^1_T=\omega_T \otimes \cI_{\Sigma}$ where $\omega_T$
is the dualizing sheaf and $\cI_{\Sigma}$ is the ideal sheaf of the nodes.

Given a bounded complex of $\cO_T$-modules
$$\cE^{\bullet}=\{0 \cdots \cE^{-p-1} \ra \cE^{-p} \ra \cE^{-p+1} \ra \cdots 0\}$$
we compute $\bR\iHom_{\cO_T}(\cE^{\bullet},\cO_T)$ using the spectral sequence
$$E_1^{p,q}=\iExt^q_{\cO_T}(\cE^{-p},\cO_T) \Rightarrow \iExt^{p+q}_{\cO_T}(\cE^{\bullet},\cO_T).$$
Note that
\begin{itemize}
\item{$\iExt^q_{\cO_T}(\phi^*\Omega^1_Y,\cO_T)=0$ for $q>0$ as $\Omega^1_Y$ is locally free;}
\item{$\iExt^q_{\cO_T}(\Omega^1_T,\cO_T)=0$ for $q>1$ by the explicit resolution above.}
\end{itemize}
In particular, only the following terms can be nonzero
$$\iHom_{\cO_T}(\phi^*\Omega^1_Y,\cO_T), \iHom_{\cO_T}(\Omega^1_T,\cO_T), \iExt^1_{\cO_T}(\Omega^1_T,\cO_T).$$
We focus on the unique interesting arrow
$$\begin{array}{rcl}
d_1: E_1^{0,0} & \ra & E_1^{1,0} \\
d\phi: \iHom_{\cO_T}(\Omega^1_T,\cO_T) & \ra & \iHom_{\cO_T}(\Omega^1_Y,\cO_T).
\end{array}
$$
Since $\iHom_{\cO_T}(\Omega^1_T,\cO_T)$ is torsion-free, $d\phi$ is injective if and only if it
is injective at generic points of $T$, which was one of our assumptions.  

Thus we have
$$E_2^{1,0}=\iHom_{\cO_T}(\phi^*\Omega^1_Y,\cO_T)/\iHom_{\cO_T}(\Omega^1_T,\cO_T)$$
and 
$$E_2^{0,1}=\iExt^1_{\cO_T}(\Omega^1_T,\cO_T).$$
Consequently 
$$\bR\iHom_{\cO_T}(\Omega_{\phi}^{\bullet},\cO_T)$$
is supported in degree one, and the associated sheaf $\cN_{\phi}$ fits into an exact sequence
$$0 \ra \iHom_{\cO_T}(\phi^*\Omega^1_Y,\cO_T)/\iHom_{\cO_T}(\Omega^1_T,\cO_T) \ra \cN_{\phi}
\ra \iExt^1_{\cO_T}(\Omega^1_T,\cO_T)\ra 0.$$
Note that the first term corresponds to deformations that leave the nodes of $T$ unchanged;
the last term is the local versal deformation space of these nodes.  
\end{proof}

\begin{rema}
In fact, $\cN_{\phi}$ is locally-free if $\phi$ is unramified (see, for example,
\cite[\S 2]{GHS})).  Conversely, if $\phi$ is ramified at a smooth point then
$\cN_{\phi}$ necessarily has torsion.  
\end{rema}

\begin{lemm} \label{lemm:rigid}
Let $\phi:T \ra S$ be an unramified morphism from a connected nodal curve of genus zero to a K3
surface.  Then 
$h^0(T,\cN_{\phi})=0$ and $h^1(T,\cN_{\phi})=1$.
In particular, the stable map $\phi:T\ra S$ is rigid in $S$.
\end{lemm}
\begin{proof}
We present an inductive argument on the number of components.  Choose a decomposition
$$T=T' \cup_p T_0, \quad T_0 \simeq \bP^1$$
where $p$ is a node disconnecting the irreducible component $T_0$ from the rest of the 
curve.  Let $\phi':T' \ra S$ be the induced map.  We have exact sequences
$$0 \ra \cN_{\phi} \otimes \cO_{T'}(-p) \ra \cN_{\phi} \ra \cN_{\phi} \otimes \cO_{T_0} \ra 0$$
and 
$$0 \ra \cN_{\phi'} \ra \cN_{\phi} \otimes \cO_{T'} \ra Q \ra 0,$$
where $Q$ is a torsion sheaf of length one supported at $p$.  
Furthermore, $\cN_{\phi} \otimes \cO_{T_0} \simeq \cO_{\bP^1}(-1)$ so the first
sequence gives
$$H^0(\cN_{\phi} \otimes \cO_{T'}(-p))=H^0(\cN_{\phi}), \quad
H^1(\cN_{\phi} \otimes \cO_{T'}(-p))=H^1(\cN_{\phi}).$$ 
The second sequence may be interpreted as $\cN_{\phi}$ tensored by
$$ 0 \ra \cO_{T'}(-p) \ra \cO_{T'} \ra \cO_p \ra 0,$$
thus we obtain $\cN_{\phi'}=\cN_{\phi} \otimes \cO_{T'}(-p)$, hence
$$H^i(\cN_{\phi'})\simeq H^i(\cN_{\phi} \otimes \cO_{T'}(-p)) = H^i(\cN_{\phi'}),
\quad i=0,1.$$
\end{proof}

\section{K3 surfaces over finite fields}
For general background and definitions, we refer the reader to \cite{RS}.

Let $S_0$ be a K3 surface over a finite field $\F_q$, and $\oS0$
the resulting surface over $\overline {\F_q}$.  Consider the
Picard group $\Pic(\oS0)$ and the $\ell$-adic
cohomology group $H^2(\oS0,\Q_{\ell}(1))$, which are related
by the cycle-class map
$$\Pic(\oS0) \ra H^2(\oS0,\Q_{\ell}(1)).$$
Frobenius acts on both these groups compatibly with this map,
and preserving the intersection form.  

The Frobenius action on $H^2(\oS0,\bQ_{\ell}(1))$ is diagonalizable over ${\overline \Q}_{\ell}$
with eigenvalues $\alpha_1,\ldots,\alpha_{22}$ \cite{DelWeilK3}.  
Since this factors through the orthogonal group, if $\alpha$ appears
as an eigenvalue then $\alpha^{-1}$ also appears.  Consequently, we
conclude that the following sets have an even number of elements 
\begin{itemize}
\item{the eigenvalues that are not roots of unity;}
\item{the eigenvalues that are roots of unity but
are not equal to $\pm 1$;}
\item{the total number of times $\pm 1$ appears as an eigenvalue.}
\end{itemize}
Using results of Nygaard and Ogus
on the Tate conjecture for K3 surfaces \cite{NyTate,NyOg} we
conclude
\begin{theo}
Let $S_0$ be a non-supersingular K3 surface over a finite field.  
Then $\Pic(\oS0)$ has even rank.
\end{theo}

Suppose $S$ is a K3 surface over a number field $F$ with integral model
$$\cS \ra \Spec(\mo_F),$$
which is smooth away from a finite set of primes.
For each finite extension $F'/F$, consider the set of primes
$$\mathrm{Ord}_{F'}(S)= \{ \mp  \in \Spec(\mo_{F'}): S_{\mp} \text{ is smooth and ordinary } \}.$$
After passing to a suitably large finite extension $F'/F$, the set
$\mathrm{Ord}_{F'}(S)$ has Dirichlet density one.  This is due to Tankeev \cite{Tankeev4}
in the special case where the Hodge group of $S_{\bC}$ is semisimple, and to
Joshi-Rajan \cite[\S 6]{JR} and Bogomolov-Zarhin \cite{BZ} in general.  

Finally, ordinary K3 surfaces are not supersingular \cite{RS2} (p. 1513 in translated version).

\section{Curves on K3 surfaces in positive and mixed characteristic}

Let $k$ be an algebraically closed field of characteristic $p>0$.
We shall require some general results on lifting to characteristic zero.  
We first review deformation and lifting results for K3 surfaces:
\begin{prop}
Let $S_0$ be a K3 surface defined over $k$.  
Then the versal deformation space of $S_0$ is smooth 
over the Witt-vectors $W(k)$ of relative dimension twenty.

Let $(S_0,f)$ be a polarized K3 surface over $k$,
such that self-intersection $f\cdot f$ is relatively prime to
the characteristic.  
Then there exists a 
polarized K3 surface $(S,f)$ over 
$W(k)$ reducing to $(S_0,f)$.  
\end{prop}
\begin{proof}
This proof is a special case of the analysis in \cite{Del}.

The deformations of $S_0$ are governed by the cohomology groups
$$H^0(S,\cT_{S_0}), \quad
H^1(S,\cT_{S_0}), \quad
H^2(S,\cT_{S_0}),$$
which parametrize infinitesimal automorphisms, infinitesimal deformations
and obstructions respectively.  
The deformation
problem for pairs $(S_0,\cL)$, where $S_0$ is a smooth projective
variety and $\cL$ is an invertible sheaf.  We have the Atiyah
extension \cite[p. 196]{IllCot} of the tangent sheaf
\begin{equation} \label{eqn:Atiyah}
0 \ra \cO_{S_0} \ra \cE_{S_0,\cL} \ra \cT_{S_0} \ra 0,
\end{equation}
classified (up to sign) by the Chern class 
$$c_1(\cL) \in H^1(S,\Omega_1^S)=\mathrm{Ext}^1(\cT_S,\cO_S),$$
which is the image of $[\cL] \in H^1(S,\cO^*_S)$ under the homomorphism
induced by 
$$d\log: \cO^*_S \ra \Omega^1_S.$$
Consider the formal deformation space of
pairs $(S,\cL)$.  The cohomology groups
$$\Gamma(S,\cE_{S,\cL}), \quad
H^1(S,\cE_{S,\cL}), \quad
H^2(S,\cE_{S,\cL})$$
parametrize infinitesimal automorphisms, infinitesimal deformations
and obstructions respectively.  

Now assume that $S$ is a K3 surface and $\cL$ is an ample
class $S$ with $[\cL]^2$ relatively prime to the characteristic. 
It follows that \cite{RS,LN,Ny}
$$\Gamma(S,\cT_S)=0.$$
Serre duality, combined with the isomorphism $\cT_S\simeq \Omega^1_S$
arising from the symplectic form, yields
$$H^2(S,\cT_S)=0.$$   Thus 
$$H^1(S,\cT_S)=\chi(S,\cT_S)=20$$
and the deformation space is smooth of this dimension over the Witt vectors. 

Furthermore, $c_1(\cL) \neq 0$ since $c_1(\cL)^2 \neq 0 \in H^2(S,\Omega^2_S)$.
Consequently, Extension~(\ref{eqn:Atiyah}) is nonsplit and the connecting
homomorphism
$$H^1(S,\cT_S) \ra H^2(S,\cO_S)$$
is nonvanishing.  Indeed, this map is just cup-product with 
$c_1(\cL)\in H^1(S,\Omega^1_S)$, which is non-zero by Serre-duality. 
It follows that
$$H^2(S,\cE_{S,\cL})=H^2(S,\cT_S)=0,$$
so deformations
of $(S,\cL)$ are unobstructed.
\end{proof}

Deligne \cite{Del} proves the following more general theorem:
\begin{theo} \label{theo:liftDeligne}
Let $(S_0,f)$ be a polarized K3 surface over an algebraically closed field $k$ of characteristic
$p$.  Then $S_0$ admits a lifting
to a possibly ramified extension of $W(k)$.  
\end{theo}
In fact, he proves that the locus $\Sigma_f$ in the formal versal
deformation space corresponding to K3 surfaces admitting $f$ as a polarization
is a Cartier divisor, not contained in the fiber over the closed point
of $\Specf(W(k))$.  
Ogus \cite[\S 2]{Ogus} has more precise lifting results for ordinary K3
surfaces.  These 
require finer analysis of Chern classes and crystalline cohomology.

\begin{coro}  \label{coro:MM}
Let $(S_0,f)$ be a polarized K3 surface over an algebraically closed field of 
characteristic $p$.  Then $S_0$ contains a rational curve.
\end{coro}
Indeed, the lifted K3 surface admits a rational curve by Theorem~\ref{theo:MM}.
This specializes to a cycle of
rational curves in characteristic $p$.

\begin{theo} \label{theo:lift}
Let $(S_0,f)$ be a polarized K3 surface over $k$.
Suppose that 
$$C=C_1+\ldots+C_r$$
is a connected union of distinct rational curves $C_i \subset S_0$,
such that $[C]$ is proportional to $f$.
Let $(S,f)$ be a polarized K3 surface over 
the Witt vectors $W(\overline k)$ reducing to $(S_0,f)$.
Assume one of the following conditions:
\begin{itemize}
\item{$S_0$ is not supersingular; or}
\item{the map from the normalization of $C$ to $S_0$ is
unramified, i.e., each branch of $C$ is nonsingular.}
\end{itemize}
Then there exists
a curve $R\subset S$, defined over a ramified finite extension of 
$W(k)$,
such that $R$ reduces to $C$ and each irreducible component of $R$ is
rational.
\end{theo}
\begin{proof} 
Consider the formal versal deformation space of $S_0$
$$\cS \ra B,$$
where $B\simeq \Specf W(k)[[x_1,\ldots,x_{20}]]$, i.e., a 
smooth formal scheme of dimension $20$ over $W(k)$.  
Let $b\in B$ denote the distinguished closed point.
For each $N\ge 1$, consider the relative stable map space
$$\ocM_{g,n}(\cS/B,Nf) \ra B.$$
This is a formal Artin stack with finite stabilizers, proper
over $B$.

We digress to explain the construction of this object.  There
are at least two possible approaches.  First, consider the category of Artinian
local rings $A$ with residue field $k$ and morphisms
$$\Sigma:=\Spec(A) \ra B.$$
Each base-change 
$$\cS_{\Sigma}:=\cS \times_{\Sigma} B \ra \Sigma$$
is proper and we can apply \cite[\S 8.4]{AV}
to show that $\ocM_{g,n}(\cS_{\Sigma}/\Sigma,Nf)$ exists.  Taking inverse
limits gives the desired formal stack.  

However, we can offer a more explicit construction.  Let $\ocM^{ps}_{g,n}$
denote the moduli stack of prestable curves of genus $g$ with $n$ marked
points, and  $\cC \ra \ocM^{ps}_{g,n}$ the universal curve 
(cf. \cite[p. 602]{Be}).  `Prestable' means nodal and connected,
but not necessarily satisfying the Deligne-Mumford stability condition;
these form an Artin stack, locally of finite type.  Standard deformation-theoretic
results (for example \cite[\S 3.2]{Sern})
show that the relative Hilbert `scheme' of a proper formal scheme
$\cY \ra B$ (parametrizing subschemes $\cZ \subset \cY$ reducing to
a closed subscheme $Z\subset \cY_b$) is prorepresentable.  
Given two proper formal schemes $\cW,\cY \ra B$,
the morphism `scheme' $\mathrm{Mor}_B(\cW,\cY)$
is also prorepresentable (cf. \cite[\S 3.4]{Sern}.)  A
stable map
$$\phi_b:(C_b,c_1(b),\ldots,c_n(b)) \ra S_b \simeq S_0$$
of genus $g$ with $n$ marked points such that $(\phi_b)_*[C_b]=Nf$,
naturally corresponds to an element of 
$$\mathrm{Mor}_B(\cC \times B,\cS) \ra \ocM^{ps}_{g,n} \times B.$$
(In general, we may have nontrivial stack structure when $\cC_b$ has 
infinite automorphisms.)
The proof of prorepresentability of
$$\ocM_{g,n}(\cS/B,Nf) \ra B$$
reduces then to the case of the morphism `scheme'.

General deformation-theoretic arguments (cf. \cite[I.2.15]{Kol}) show that
the relative dimension of $\ocM^{\circ}_{0,0}(\cS/B,Nf)$
over $B$ at $\phi$ is at least
$$\chi(T,\cN_{\phi})+\dim(B)=\dim(B)-1.$$
When $\phi:T \ra S_0$ is unramified,
Lemma~\ref{lemm:rigid} implies
that the fibers of $\ocM^{\circ}_{0,0}(\cS/B,Nf) \ra B$ are
zero-dimensional.  Otherwise, we use the hypothesis that
$S$ is not supersingular and $\phi$ is {\em generically} unramified to
conclude that $\phi$ does not deform to another genus-zero stable map. 
In either case, the dimensions of $\ocM^{\circ}_{0,0}(\cS/B,Nf)$ and its
image in $B$ are at least $20$.  
On the other hand, this image is contained in the locus 
$$\Sigma_{Nf} \subset B$$
parametrizing K3 surfaces admitting $Nf$ as a polarization. 
Indeed, in each fiber
$$(\phi_t)_*\cC_t =Nf.$$

By Theorem~\ref{theo:liftDeligne}, the
formal scheme $\Sigma_{Nf}$  has dimension $20$
and is not contained in the fiber over the closed point of $\Specf(W(k))$.
The same must hold for $\ocM^{\circ}_{0,0}(\cS/B,Nf)$, so there
are formal lifts of $\phi:T\ra S$ to genus-zero maps in characteristic zero.

It remains to show these formal deformations are algebraic.  For this
purpose, we restrict to the polarized deformation space
$$\cS_{\Sigma_{Nf}} \ra \Sigma_{Nf},$$
which is projective in the sense that it admits a formal embedding
into a projective space $\bP^d_{\Sigma_{Nf}}, d=\chi(S,\cO_S(Nf))-1.$
This deformation is algebraizable by standard results of Grothendieck
(see \cite[2.5.13]{Sern}, for example.)  It follows that the
associated moduli spaces of stable maps is algebraizable as well.  Indeed,
moduli spaces of stable maps into projective schemes are proper
stacks with projective coarse moduli spaces.  
\end{proof}

\section{Proof of the Main Theorem}
Assume $\Pic(S)$ is generated by an ample class $f$, of arbitrary degree.
Suppose $S$ admits a finite number of rational curves $R_1,\ldots,R_s$
with classes $[R_i]=m_if$ and write $m=\max\{m_1,\ldots,m_s \}$.  

\begin{lemm} \label{lemm:help}
There are only a finite number of primes $p$ such that there exists 
a curve $C \subset S_p$ with $[C] \not \in \bZ f$ and 
$$C\cdot f \le m f\cdot f.$$
\end{lemm}
\begin{proof}
There are a finite number of rank-two extensions 
$$\bZ f \subset \Lambda,$$
where $\Lambda$ is an integral lattice of signature $(1,1)$ admitting
a vector $v$ linearly independent from $f$ with $v\cdot f \le m f\cdot f$.
For each such lattice, there are at most a finite number of primes $p$
such that $\Lambda \subset \Pic(S_p)$.  
\end{proof}

Now we assume $S$ has degree two.
Let $\iota:S \ra S$ denote the involution associated to the 
branched double cover $S \ra \P^2$.  It acts on the primitive
cohomology of $S$ via multiplication by $-1$.  Let $\iota_p$ denote
its reduction mod $p$.
We shall derive a contradiction by producing an irreducible
rational curve in a class
$Nf$ for some $N>m$.  

Choose $p$ to a prime of good reduction that is not in the finite
set of primes specified in Lemma~\ref{lemm:help}.  Let $N$ be the 
the smallest positive integer such that $Nf$ is decomposable
$$Nf=a_1[C_1]+\ldots+a_r[C_r],$$
where the $a_i$ are positive integers and not all of the $[C_i]$ are
proportional to $f$.   We assume each class $[C_i]$ is indecomposable
and the $C_i$ is irreducible and rational.  (This is possible by
the Mori-Mukai argument.)  By the minimality of $N$, none of the $C_i$
is proportional to $f$.  

Since
$$\iota_p^*C_i+C_i=C'_i + C_i$$
is invariant under $\iota_p$ it equals some multiple of $f$.  We claim
this lifts to an irreducible rational curve $R$ over $\overline{\Q}$.  
Consider the chain of two $\bP^1$'s
$$T=\{xy=0 \} \subset \bP^2$$
and choose a birational morphism $\phi:W \ra C'_i+C_i$.  
Let $j:T\ra S_p$ be the induced morphism;  then $T$ is the 
specialization of a rational curve over the Witt-vectors by Theorem~\ref{theo:lift}.
Since this curve is rigid, it can be defined over $\overline{\Q}$ as well.

\bibliographystyle{smfplain}
\bibliography{rank1}

\end{document}